\documentclass[12pt]{amsart}

\usepackage{latexsym}
\usepackage{graphics}
\usepackage{amsmath,amssymb,amscd}
\usepackage[all]{xy}

\newtheorem{proposition}{Proposition}[section]
\newtheorem{theorem}[proposition]{Theorem}
\newtheorem{lemma}[proposition]{Lemma}

\newcommand{\fq}{\ensuremath{\mathbb{F}_q}}
\newcommand{\lra}{\longrightarrow}

\newcommand{\Map}{\mathop{\mathrm{Map}}}
\newcommand{\Hom}{\ensuremath{\mathrm{Hom}}}
\newcommand{\A}{\ensuremath{\mathcal{A}}}

\newcommand{\ddr}[5]{\ensuremath{#1\stackrel{#2}
{\longrightarrow}#3\stackrel{#4}{\longrightarrow}#5}}

\newcommand{\Spin}{\ensuremath{\mathrm{Spin_7}}}
\newcommand{\Sol}{\ensuremath{\mathrm{Sol}}}
\newcommand{\DI}{\ensuremath{\mathrm{DI}}}
\newcommand{\G}{\ensuremath{\mathrm{G}}}

\begin{document}
\title{The cohomology of exotic \\$2$--local finite groups}

\author{Jelena Grbi\' c}


\maketitle
\begin{abstract}
There exist spaces $B\Sol(q)$ which are the classifying spaces of a
family of 2--local finite groups based on certain fusion system over
the Sylow 2--subgroups of $\Spin(q)$. In this paper we calculate the
cohomology of $B\Sol(q)$ as an algebra over the Steenrod algebra
$\mathcal{A}_2$. We also provide the calculation of the cohomology
algebra over $\mathcal{A}_2$ of the finite group of Lie type
$\G_2(q)$.
\end{abstract}

\section{Introduction}
\label{intro}

In 1974, working on the classification of finite groups, Solomon
\cite{So} looked at {\it fusion systems} over a finite group $S$,
that is, a category whose objects are the subgroups of $S$ and whose
morphisms are monomorphisms of groups which include all those
induced by conjugation by elements of $S$. In particular, he
considered groups in which all involutions are conjugate, and such
that the centralizer of each involution contains a normal subgroup
isomorphic to $\Spin(q)$ with odd index, where $q$ is an odd prime
power. Solomon showed that a group with such a fusion system does
not exist. Later on, in 1996 Benson~\cite{Be}, inspired by Solomon's
work, constructed a family of spaces that he informally entitled
{\it``$2$--completed classifying spaces of the nonexistent Solomon
finite groups"}, and denoted them by $B\Sol(q)$. Chronologically
first Dwyer and Wilkerson \cite{DW} came up with the celebrated
discovery of the exotic finite loop space $\DI(4)$. They constructed
a space $B\DI(4)$ which realises the mod~$2$ Dickson invariants of
rank 4, that is,
\[
H^*(B\DI(4);\mathbb{F}_2) \cong P[u_8,u_{12},u_{14},u_{15}],
\]
with the action of the Steenrod algebra $\A_2$ given by
$Sq^4u_8=u_{12}$, $Sq^2u_{12} = u_{14}$ and $Sq^1u_{14}=u_{15}$.
Soon after the experts became aware of the existence of self
homotopy equivalences $\Psi^q\colon B\DI(4)\lra B\DI(4)$ for each
odd prime power~$q$, which are analogues of the unstable Adams
operations on classifying spaces of compact Lie groups completed
away from $q$. The first publish proof of this fact appeared in
Notbohm's work \cite{No}. Benson observed that Dwyer-Wilkerson's
finite loop space $\DI(4)$ carries a $2$--local structure closely
related to Solomon's fusion systems. For each odd prime power $q$,
Benson~\cite{Be} defined a space $B\Sol(q)$ as the homotopy pullback
space of the diagram
\begin{equation}\label{d1}
\xymatrix{B\Sol(q) \ar[r] \ar[d] & B\DI(4) \ar[d]^{\Delta} \\
           B\DI(4)\ar[r]^-{(1, \Psi^{q})} & B\DI(4)\times B\DI(4).}
\end{equation}

Benson claimed that the space $B\Sol(q)$ realises the $2$--fusion
pattern over the Sylow 2--subgroup of $\Spin(q)$ studied by
Solomon~\cite{So} and shown by him not to correspond to any genuine
finite group.

In~\cite{LO} Levi and Oliver have been interested in the fusion
patterns studied by Solomon in the context of saturated fusion
systems and ${p\text{--local}}$ finite groups. They have shown that
the fusion system studied by Solomon is saturated in the sense of
Puig~\cite{Pu} and that each one of these fusion systems gives rise
to a unique 2--local finite group, which they named $\Sol(q)$. Their
approach to the subject enabled them to show that if $n|m$, then
$\Sol(q^n)$ can be thought of as a subgroup of $\Sol(q^m)$ and that
the resulting limiting term, which they denoted by $\Sol(q^\infty)$,
has a classifying space homotopy equivalent to $B\DI(4)$. This in
turn made it possible to give an alternative definition of the
unstable Adams operations $\Psi^q$ on $B\DI(4)$ and to show that the
space $B\Sol(q)$ constructed by Benson is homotopy equivalent to the
classifying space of the $2$--local finite group $\Sol(q)$, for each
odd prime power $q$.

The object of this paper is the calculation of the cohomology
algebra over the Steenrod algebra $\mathcal{A}_2$ of ``Solomon's
groups". We use Benson's notation, denoting the classifying space
of Solomon's groups by $B\Sol(q)$. The cohomology of a $p$--local
finite group \cite{BLO} is given in the appropriate sense as the
subalgebra of ``stable elements'' in the cohomology of its Sylow
$p$--subgroup. However, as in the case of ordinary finite groups,
calculating cohomology by stable elements is quite a difficult
task. The fact that $B\Sol(q)$ is given by a rather simple
pullback diagram makes the calculation of the cohomology as a
module over the mod $2$ Dickson invariants of rank $4$, done
initially by Benson \cite{Be}, straightforward. The difficulties
arise when one wants to determine the algebra structure of the
cohomology of $B\Sol(q)$ with $\mathbb{F}_2$ coefficients as an
algebra over $\mathcal{A}_2$.

The main result of the paper is given by the following theorem.
\begin{theorem}
\label{Thm1}
 Fix an odd prime power $q$. Then
\[
H^{\ast}(B\Sol(q); \mathbb{F}_{2})\cong \mathbb{F}_{2}[u_{8},u_{12},
u_{14}, u_{15}, y_{7},y_{11},y_{13}]/I,
\]
where $I$ is the ideal generated by the polynomials
\begin{enumerate}
 \item[] $y_{11}^{2}+u_{8}y_{7}^{2}+u_{15}y_{7}$,
 \item[] $y_{13}^{2}+u_{12}y_{7}^{2}+u_{15}y_{11}$,
 \item[] $y_{7}^{4}+u_{14}y_{7}^{2}+u_{15}y_{13}$.
\end{enumerate}
The action of the Steenrod algebra $\mathcal{A}_2$ is determined by
\[
Sq^{4}(u_{8})=u_{12}, \quad Sq^1(u_{12})=0,\quad
Sq^{2}(u_{12})=u_{14},\quad Sq^{1}(u_{14})=u_{15},
\]
\[
Sq^1(y_7)=0, \quad Sq^{4}(y_{7})=y_{11}, \quad Sq^1(y_{11})=0,\quad
Sq^{2}(y_{11})=y_{13},\quad Sq^{1}(y_{13})=y_{7}^{2}
\]
and the Steenord algebra axioms.

The second page of the Bockstein spectral sequence is given by
\[
B_2(B\Sol(q))\cong P[u_8,u_{12},u_{14}^2]\otimes
E[y_7,y_{11},y_7^2y_{13}].
\]
The higher Bockstein operators are given as follows. Let
$k=\nu_2(q^2-1)$, which is at least $3$. Then
\[
\beta_{k+1}(y_7)=u_8,\quad \beta_k(y_{11})=u_{12},\quad
\beta_k(y_7u_8u_{12})=u_8^2u_{12}, \quad
\beta_{k-1}(y_7^2y_{13}+\epsilon u_8^2y_{11})= u_{14}^2
\]
where $\epsilon=0$ or $1$.
\end{theorem}

As regards the cohomology of the classifying space of the finite
group of Lie type $\G_2(q)$ with $\mathbb{F}_2$ coefficients,
Kleinerman \cite{Kl} computed its module structure, while Milgram
\cite{Mi} described it as an algebra over the Steenrod algebra.
Although Milgram's work is elementary in a sense, it is quite
lengthy and technical. The method that has been used to calculate
the cohomology algebra of $B\Sol(q)$ can be also applied to the
calculation of the cohomology algebra over the Steenrod algebra
$\mathcal{A}_2$ of $B\G_2(q)$. As it is based in deeper results
(Friedlander's fibre square) it produces a short and elegant proof
that we shall present in the Appendix. The result is as follows.
\begin{theorem}
Fix an odd prime power $q$. Then
\[
H^*(B\G_2(q); \mathbb{F}_2)\cong \mathbb{F}_2[d_4,d_6,d_7,y_3,y_5]/I
\]
where $I$ is the ideal generated by the polynomials
\begin{enumerate}
 \item[] $y_5^2+y_3d_7+y_3^2d_4$
 \item[] $y_3^4+y_5d_7+y_3^2d_6$.
\end{enumerate}
The action of $\A_2$ is determined by
\[
Sq^1d_4=0,\, Sq^2d_4=d_6,\, Sq^1d_6=d_7,\, Sq^1y_3=0,\,
Sq^2y_3=y_5,\, Sq^1y_5=u_3^2
\]
and the Steenord algebra axioms.
\end{theorem}
For the sake of completeness, we state Kleinerman's result \cite{Kl}
on the higher Bockstein relations in $H^*(B\G_2(q))$. The second
page of the Bockstein spectral sequence is given by
\[
B_2(B\G_2(q))\cong P[d_4,d_6^2]\otimes E[y_3, y_3y_5].
\]
The higher Bockstein operators are given as follows. Let
$k=\nu_2(q^2-1)$. Then
\[
\beta_k(y_3)=d_4\quad \text{ and }\quad
\beta_{k-1}(y_3^2y_5)=d^2_6+\epsilon d_4^3
\]
for $\epsilon=0$ or $1$.

Through the rest of the paper, if not specified differently, we work
over the field $\mathbb{F}_2$ and denote by $H^*(X)$ the cohomology
of a topological space $X$ with $\mathbb{F}_2$ coefficients.

\noindent\textbf{Acknowledgements.} The author would like to thank
the referee for many important and useful suggestions and comments.
The manuscript has been substantially modified, in particular
section 4, following the referee's suggestions.

\section{The Cohomology of $B\Sol(q)$ as a Module over
$\mathcal{A}_2$} \label{sec:1}

In this section we apply the Eilenberg--Moore spectral sequence to
pullback diagram \eqref{d1} to derive the cohomology of $B\Sol(q)$
as a module over the mod~$2$ Steenrod algebra $\mathcal{A}_2$.

Since the homotopy fibre of the diagonal map $\Delta\colon X\lra
X\times X$ is $\Omega X$, and $B\Sol(q)$ is defined by homotopy
pullback diagram \eqref{d1}, the homotopy fibre of $p\colon
B\Sol(q)\lra B\DI(4)$ is homotopy equivalent to $DI(4)$ and there is
a map of homotopy fibrations
\begin{equation}
\label{maindgrm}
\xymatrix{\DI(4)\ar[d]\ar@{=}[r] & \DI(4)\ar[d]\\
B\Sol(q) \ar[r] \ar[d]^p & B\DI(4) \ar[d]^{\Delta} \\
           B\DI(4)\ar[r]^-{(1, \Psi^{q})} & B\DI(4)\times B\DI(4)}
\end{equation}
where $\DI(4)$ is the loop space $\Omega B\DI(4)$. The following
proposition states the $H^*(B\DI(4))$--module structure of
$H^*(B\Sol(q))$.
\begin{proposition}
\label{propmodule} The mod~$2$ cohomology of $B\Sol(q)$ as a
$P[u_{8},u_{12}, u_{14}, u_{15}]$--module is given as
\[
H^*(B\Sol(q);\mathbb{F}_2)\cong P[u_{8},u_{12}, u_{14},
u_{15}]\otimes E[v_{11},v_{13}]\otimes P[v_7]/v_7^4.
\]
The action of the Steenrod algebra $\mathcal{A}_2$ is determined
by
\[
Sq^{4}(u_{8})=u_{12}, \quad Sq^1(u_{12})=0,\quad
Sq^{2}(u_{12})=u_{14},\quad Sq^{1}(u_{14})=u_{15},
\]
\[
Sq^1(v_7)=0, \quad Sq^{4}(v_{7})=v_{11}, \quad
Sq^1(v_{11})=0,\quad Sq^{2}(v_{11})=v_{13},\quad
Sq^{1}(v_{13})=v_{7}^{2}
\]
and the Steenord algebra axioms.
\end{proposition}

\begin{proof}
We look at the Eilenberg--Moore spectral sequence of pullback
diagram~\eqref{d1}. The above spectral sequence satisfies the
hypothesis of Smith's ``big collapse theorem''
\cite[Theorem~II.3.1]{Sm} and hence collapses at the $E_2$~term.
Namely, one obtains
\[
E_\infty\cong E_2\cong P[u_8, u_{12}, u_{14}, u_{15}]\otimes
E[a_7, a_{11}, a_{13}, a_{14}]
\]
with $P[u_8, u_{12}, u_{14}, u_{15}]\cong H^*(B\DI(4))$ in
filtration degree $0$ and $a_i$'s in filtration degree $-1$ with the
Steenord algebra action $Sq^4(a_7)=a_{11}, Sq^2(a_{11})=a_{13}$ and
$Sq^1(a_{13})=a_{14}$. If we take $v_7\in H^7(B\Sol(q))$
representing $a_7$ and then $v_{11}=Sq^4(a_7), v_{13}=Sq^2(v_{11})$
and $v_{14}=Sq^1(v_{13})$, then the Adem relations show that
$v_{14}=v_7^2$. That proves the proposition.
\end{proof}

We want to make the statement of the last proposition more
geometrical. Namely, we want to show that
\[
H^*(B\Sol(q))
\cong H^*(B\DI(4)) \otimes H^*(\DI(4)),
\]
as $H^*(B\DI(4))$--modules and that the action of the Steenord
algebra splits.

By the following proposition, that will have a crucial role later on
in determining the algebra structure of $H^*(B\Sol(q))$, we
investigate the interaction of the action of the Steenrod algebra on
$H^*(B\DI(4))$ and $H^*(\DI(4))$ in the category of modules over
$\mathcal{A}_2$.

\begin{proposition}
\label{Lannes} The factor $H^*(B\DI(4))$ is a split summand of
$H^*(B\Sol(q))$ in the sense that the natural map $H^*(B\DI(4))\lra
H^*(B\Sol(q))$ has a left inverse in the category of unstable
algebras over $\A_2$.
\end{proposition}

\begin{proof}
From the construction of $\DI(4)$ in \cite{DW}, there is a unique,
up to homotopy, map $f\colon B(\mathbb{Z}/2)^4\lra B\DI(4)$ such
that $f^*\colon H^*(B\DI(4))\lra H^*(B(\mathbb{Z}/2)^4)$ as a map of
$\A_2$ modules, is the inclusion of mod~2 Dickson invariants of
rank~4. As mentioned before, $\Psi ^q$ is the identity in mod~2
cohomology, so $\Psi ^q\circ f\simeq f$ by the uniqueness property
of $f$. Thus $(1,\Psi ^q)\circ f$ lifts through $\Delta$, that is,
$(1,\Psi ^q)\circ f\simeq \Delta\circ f$. Now, using the universal
property of pullback diagram \eqref{maindgrm}, there is a map
$g\colon B(\mathbb{Z}/2)^4 \lra B\Sol(q)$ such that $p\circ g\simeq
f$. Thus in cohomology the images of both maps $f^*$ and $g^*\circ
p^*$ are the mod~$2$ Dickson invariants of rank~$4$ sitting inside
$H^*(B(\mathbb{Z}/2)^4)$. Therefore the natural map
${H^*(BDI(4))\longrightarrow H^*(B\Sol(q))}$ has a left inverse in
the category of unstable algebras over~$\mathcal{A}_2$.
\end{proof}

Now consider the path-loop fibration sequence
\[
\ddr{\DI(4)}{}{\mathcal{P}B\DI(4)}{}{B\DI(4)}.
\]
Applying the Eilenberg--Moore spectral sequence to this fibration,
we obtain the module structure of $H^*(\DI(4))$ that is isomorphic
to $E[x_7,x_{11},x_{13},x_{14}]$. Since in this case the
Eilenberg--Moore spectral sequence is compatible with $\A_2$
\cite{Sm70}, one can also determine the action of $\A_2$ as
$Sq^4x_7=x_{11}$, $Sq^2x_{11}= x_{13}$ and $Sq^1x_{13}=x_{14}$. A
Steenrod squares manipulation now implies that
\[
H^*(\DI(4))\cong P[x_7]/x_7^4 \otimes E[x_{11},x_{13}],
\]
as an algebra.

This shows that \[ H^*(B\Sol(q)) \cong H^*(B\DI(4)) \otimes
H^*(\DI(4)),
\]
as an $H^*(B\DI(4))$--module; that the action of the Steenord
algebra splits and that the cohomology classes $v_i$ in
$H^*(B\Sol(q))$ are detected in $H^*(\DI(4))$.
\section{The Cohomology of $B\Sol(q)$ as an Algebra over
$\mathcal{A}_2$}
 \label{sec:2}

This section is devoted to the calculation of the algebra structure
of $H^*(B\Sol(q))$. By the previous section,
\[
H^*(B\Sol(q)) \cong P[u_8,u_{12},u_{14},u_{15}]\otimes
P[v_7]/(v^4_7)\otimes E[v_{11},v_{13}]
\]
as modules over $P[u_8,u_{12},u_{14},u_{15}]$. Thus, let $y_7\in
H^*(B\Sol(q))$ denote the unique generator which maps to $v_7\in
H^*(\DI(4))$ under the map induced by the fibre inclusion. Let
$y_{11}, y_{13}\in H^*(B\Sol(q))$ denote $Sq^4(y_7)$ and
$Sq^2Sq^4(y_7)$ respectively. The structure given above is also
the algebra structure of $H^*(B\Sol(q))$ if and only if the
subalgebra of $H^*(B\Sol(q))$ generated by $y_7, y_{11}$ and
$y_{13}$ is isomorphic to $P[v_7]/(v^4_7)\otimes
E[v_{11},v_{13}]$. As we shall show, whether or not this is the
case depends only on a single parameter $A\in\mathbb{F}_2$, the
value of which determines all the relations in $H^*(B\Sol(q))$. We
shall then proceed by showing that the cohomology algebra cannot
split as a tensor product, thus determining the value of $A$ to be
1. The full structure of $H^*(B\Sol(q))$ as an algebra over the
Steenrod algebra will then follow at once.

Consider the possible multiplications on $H^*(B\Sol(q))$. Looking at
the action of the Steenrod algebra $\A_2$, we have $ y_{11}^2 =
Sq^8(y_7^2)$. Because of dimensional reasons, $ y_{11}^2$ can be
presented as
\begin{equation}
\label{rel1}
 y_{11}^2 = Ay_7^2u_8 + By_7u_{15} +
Cu_8u_{14}.
\end{equation}
On the one hand we have $Sq^1 y_{11}^2=0$, and on the other hand
applying $Sq^1$ to relation \eqref{rel1}, we have $Sq^1y_{11}^2 =
Cu_8u_{15}$. This implies that $C=0$. Next
\[
y_{13}^2 = Sq^4y_{11}^2 = Ay_7^2u_{12} + By_{11}u_{15}.
\]
Notice that $Sq^1y_{13}=y^2_7$ as $y_{13}=Sq^2Sq^4y_7$. Hence
\[
y_7^4 = Sq^2y_{13}^2 = Ay_7^2u_{14} + By_{13}u_{15}
\]
and
\[
0=Sq^1y_7^4 = (A+B)y_7^2u_{15}.
\]
This implies that $A=B$ and so one derives the relations claimed in
the main theorem if $A=1$. Notice also that if $A=0$, then there is
a ring homomorphism $H^*(\DI(4))\lra H^*(B\Sol(q))$ and the
composite
\begin{equation}
\label{splitting} H^*(B\DI(4))\otimes H^*(\DI(4))\lra
H^*(B\Sol(q))^{\otimes 2} \lra H^*(B\Sol(q))
\end{equation}
is an isomorphism of vector spaces and a map of algebras. Hence
composite~\eqref{splitting} is an isomorphism of algebras over the
Steenrod algebra. Thus to prove the first part of Theorem $1$ (not
including the information on the Bockstein spectral sequence) it
remains to show the following lemma.
\begin{lemma}
\label{product} The cohomology  $H^*(B\Sol(q))$ cannot split as a
tensor product of $H^*(B\DI(4))$ and $H^*(\DI(4))$ in the category
of algebras over~$\A_2$.
\end{lemma}

\begin{proof}
Let $V=(\mathbb{Z}/2)^4$, $f\colon BV\lra B\DI(4)$ the maximal
elementary abelian $2$--subgroup of $\DI(4)$ and $g\colon BV\lra
B\Sol(q)$ the lift constructed in the proof of
Proposition~\ref{Lannes}.

Applying the functor $\Map(BV,-)$ to the homotopy pullback diagram
\[
\xymatrix{B\Sol(q) \ar[r] \ar[d] & BDI(4) \ar[d]^{\Delta} \\
           BDI(4)\ar[r]^-{(1, \Psi^{q})} & BDI(4)\times BDI(4)}
\]
and fixing the component of $f$, we obtain the homotopy pullback
diagram
\begin{equation}
\label{mapdgm} \xymatrix{
  \Map(BV,B\Sol(q))_F\ar[r] \ar[d] & \Map(BV,B\DI(4))_f \ar[d]^{\Delta}  \\
  \Map(BV, B\DI(4))_f \ar[r]^-{\Delta} & \Map(BV,B\DI(4))_f^{\times 2}}
\end{equation}
where $F$ stands for all of the connected components of all lifts
(up to homotopy) of $f:BV\lra B\DI(4)$ to $B\Sol(q)$. It is shown in
\cite{LO} that there are more than just one faithful representation
of an elementary abelian $2$--group of rank $4$ in Solomon
$2$--local finite group.

Let us assume that $H^*(B\Sol(q))\cong H^*(B\DI(4))\otimes
H^*(\DI(4))$ as algebras over the Steenord algebra. Then there is
only one lift $g$ of $f$ contradicting the existence of pullback
diagram \eqref{mapdgm}. This is based on the fact that
\begin{equation*}
\begin{split}
[BV, B\Sol(q)]&\cong \Hom_\A(H^*(B\Sol(q)),H^*(BV)) \\
              &\cong \Hom_\A(H^*(B\DI(4)))\otimes H^*(\DI(4)),
              H^*(BV))\\
              &\cong \Hom_\A(H^*(B\DI(4)),H^*(BV))
\end{split}
\end{equation*}
where the first isomorphism is given in \cite{La}, and the last
isomorphism follows since $H^*(\DI(4))$ is a finite vector space.
Further on, Dwyer and Wilkerson \cite{DW} showed that
$\Hom_\A(H^*(B\DI(4)),H^*(BV))$ contains only one non-identity
homotopy class. That proves the assertion that there is at most one
lifting $g$ of $f$ and finishes the proof of the lemma.
\end{proof}

\section{The Bockstein Spectral Sequence for $B\Sol(q)$}
\label{sec:3} In this section we describe the Bockstein operators of
$H^*(B\Sol(q))$.

As it has been seen, only the module structure is needed for the
$Sq^1$ calculation as the algebra splits the same way as a module
over $Sq^1$ (Proposition~\ref{propmodule}). Thus the second page of
the Bockstein spectral sequence looks like
$$
B_2(B\Sol(q))\cong P[u_8,u_{12},u_{14}^2]\otimes
E[v_7,v_{11},v_7^2v_{13}].
$$
A Serre spectral sequence calculation in cohomology with
coefficients in the $2$--adic integers shows that $v_7$, $v_{11}$,
$u_8u_{12}v_7$ and $v_7^2v_{13}$ support higher Bocksteins.

Notbohm \cite{No} showed that the effect of the Adams operation
$\Psi^q$ from diagram \eqref{maindgrm} on
$H^{2n}(B\DI(4);\hat{\mathbb{Z}}_2)\otimes\mathbb{Q}$ is
multiplication by $q^n$. By \cite[Proposition~3.3]{No}, the Adams
operation $\Psi^q$ is an equivalence and therefore it induces an
isomorphism of the algebra $H^*(B\DI(4); \mathbb{F}_2)$. For
dimensional reasons $(\Psi^q)^*(u_8)=u_8$, so $(\Psi^q)^*$ is the
identity on mod~$2$ cohomology. Now, apply the Serre spectral
sequence with $2$--adic coefficients to the fibration
\begin{equation}
\label{firstfib} \ddr{\DI(4)}{}{B\DI(4)}{\Delta}{B\DI(4)\times
B\DI(4)}.
\end{equation}
Notice first that $\Delta^*$ is epimorphic, thus all classes above
the bottom line of the spectral sequence are annihilated by a
differential, that is, $E^{p,q}_\infty=0$ for $q>0$. Furthermore the
classes $E^{p,0}_2$ that are in the image of some differential are
those that belong to the kernel of $\Delta^*$. Let $\rho\colon
\hat{\mathbb{Z}}_2 \lra \mathbb{F}_2$ be the reduction mod $2$ map.
Denote by $y_i$ classes in $H^*(\DI(4);\hat{\mathbb{Z}}_2)$ such
that $\rho(y_i)=v_i$ and by $\omega_j$ classes in
$H^*(B\DI(4);\hat{\mathbb{Z}}_2)$ such that $\rho(\omega_j)=u_j$.
For degree reasons $y_7\in H^*(DI(4); \hat{\mathbb{Z}}_2)$
transgresses to $\omega_8-\omega'_8\in
H^*(B\DI(4);\hat{\mathbb{Z}}_2)^{\otimes 2}$. Further on using
pullback diagram~\eqref{maindgrm}, the class $\omega_8-\omega'_8$,
under the effect of the map $(1,\Psi^q)^*$, maps on
$(q^4-1)\omega_8\in H^*(B\DI(4);\hat{\mathbb{Z}}_2)$. Now using the
naturality of the Serre spectral sequence, in the homotopy fibration
\begin{equation}
\label{secondfib}
\ddr{\DI(4)}{}{B\Sol(q)}{}{B\DI(4)}
\end{equation}
the integral cohomology class $y_7\in
H^*(\DI(4);\hat{\mathbb{Z}}_2)$ transgresses onto the class
$(q^4-1)\omega_8\in H^*(B\DI(4);\hat{\mathbb{Z}}_2)$.

In an analogue way, the class $y_{11}$ transgresses onto the class
$(q^6-1)\omega_{12}\in H^*(B\DI(4);\hat{\mathbb{Z}}_2)$.

For degree reasons, in the Serre spectral sequence for
fibration~\eqref{firstfib} the first non-trivial differential on the
class $y_7\omega_8\omega_{12}\in E_2^{20,7}$ is $d_8$ so that
\[
d_8(y_7\omega_8\omega_{12})=d_8(y_7)\omega_8\omega_{12}=\omega_8^2\omega_{12}.
\]
Further on using pullback diagram~\eqref{maindgrm}, the class
$\omega_8^2\omega_{12}$ under the effect of the map $(1,\Psi^q)^*$,
maps on $(q^{14}-1)\omega_8^2\omega_{12}\in
H^*(B\DI(4);\hat{\mathbb{Z}}_2)$. Therefore in the Bockstein
spectral sequence for $B\Sol(q)$,
$\beta_{\nu_2(q^{14}-1)}(v_7u_8u_{12})=u_8^2u_{12}$.

Now we want to know what happens with the integral class
$v_7^2v_{13}$ of $H^*(\DI(4))$. Denote by $y_{27}$ the class in
$H^*(\DI(4); \hat{\mathbb{Z}}_2)$ such that
$\rho(y_{27})=v_7^2v_{13}$. In the Serre spectral sequence of
fibration \eqref{firstfib}, the first possibly non-trivial
differential on $y_{27}$ is $d_{14}$. In particular,
\[
\begin{array}{ll}
E_2^{14,14}& =H^{14}(B\DI(4)\times B\DI(4);
H^{14}(\DI(4))) \\
   & =H^{14}(B\DI(4)\times B\DI(4);
\mathbb{Z}/2)=\mathbb{Z}/2\oplus\mathbb{Z}/2
\end{array}
\]
having two $2$--torsion classes $y_7^2(\omega_{14}\otimes 1)$ and
$y_7^2(1\otimes \omega_{14})$. These classes can be only hit by a
differential coming from $E_{15}^{0,27}$. Since
$E^{0,27}_{15}=y_{27}\hat{\mathbb{Z}}_2$ and $E^{29,0}_{16}$
contains only the extension class
$\omega_{29}\mathbb{Z}_2=\mathrm{Ext}(\omega_{14}\mathbb{Z}_2,
\omega_{14}\mathbb{Z}_2)$ it follows that $d_{14}(y_{27})$ is
nontrivial and $d_{15}\colon E_{15}^{14,14}\lra E_{15}^{29,0}$ is an
isomorphism. The above argument shows that in page 28 we are left
with $2y_{27}\hat{\mathbb{Z}}_2$ in $E^{0,27}_{28}$, hence
$d_{28}(2y_{27})=\omega_{14}^2\otimes 1-1\otimes
\omega_{14}^2$.

Using the information obtained in the above paragraphs, we write the
Serre spectral sequence for $\ddr{\DI(4)}{}{B\Sol(q)}{}{B\DI(4)}$
with $2$--adic coefficients. The pullback diagram \eqref{maindgrm}
provides a maps of Serre spectral sequences. First notice that there
is a $2$--torsion class in $E^{14,14}_2$ which must survive to
$E_\infty$ as it is the only class in total degree 28, and according
to the calculation of the Bockstein spectral sequence, there is a
$2$--torsion class in $H^{28}(B\Sol(q);\hat{\mathbb{Z}}_2$. Thus all
previous differentials are trivial on $y_{27}$, and therefore
\[
d_{28}(y_{27})=\frac{q^{14}-1}{2}\omega^2_{14}.
\]
So there must be
$\overline{v_{27}}\in H^*(B\Sol(q);\mathbb{Z}/2)$ with
$\beta_{\nu_2(q^{14}-1)-1}(\overline{v_{27}})=\omega_{14}^2$. It
follows that $\overline{v_{27}}=v_7^2v_{13}+\epsilon u_8^2v_{11}$,
where $\epsilon=0$ or $1$.

The following remark is kindly suggested by the referee in order to
underline the variety of approaches that could be taken to solve the
above problem.

\noindent \textbf{Remark.} Since we have the transgression
$d_{28}(y_{27})=\frac{q^{14}-1}{2}\omega^2_{14}$, the method of
universal example (construct a fibration of Eilenberg-MacLane spaces
detecting the relevant cohomology classes), shows that
$\beta_{(\nu_2(q^{14}-1)-1)}(y_7^2v_{13})=\omega_{14}^2$. In fact,
also $\beta_{(\nu_2(q^{14}-1)-1)}(y_7^2v_{13}+
\omega_8^2y_{11})=\omega_{14}^2$.

We finish the proof by showing that in $H^*(B\Sol(q))$ for every
$k$, ${d_k(v_7^2)=0}$. First of all, direct calculation shows that
$d_8(y_7^2)=0$ because $d_8$ is an antiderivation. Since
$\omega_{15}$ is a $2$--torsion element and $\Psi^q$ is an
isomorphism, it follows that $\Psi^q(\omega_{15})=\omega_{15}$. As
in the Serre spectral sequence of fibration~\eqref{secondfib} with
coefficients in $\mathbb{F}_2$ the class $\omega_{15}$ survives, we
have $d_{15}(v^2_7)=0$.

Therefore we have determined that in the mod~$2$ Bockstein
spectral sequence the pairs $(v_7,u_8)$, $(v_{11},u_{12})$,
$(v_7u_8u_{12},u_8^2u_{12})$ and $(v_7^2v_{13}+\epsilon
u_8^2v_{11}, u_{14}^2)$ are connected by higher Bocksteins
operations of orders $\nu_{2}(q^{4}-1),
\nu_{2}(q^{6}-1),\nu_{2}(q^{14}-1)$ and $\nu_{2}(q^{14}-1)$,
respectively.

\section{Appendix}
\label{sec:4} In this Appendix we give the cohomology of the
finite group of Lie type $\G_2(q)$ as an algebra over the Steenrod
algebra $\A_2$ using the methods developed for the cohomology
calculation of $Sol(q)$. The module structure of the cohomology of
$B\G_2(q)$ is well-known (see for example \cite{Kl}). As already
mentioned in the introduction Milgram calculated the algebra
structure of $H^*(B\G_2(q))$ using completely different methods to
ours.

\begin{theorem}
\label{cohg2q} Fix an odd prime power $q$. Then
\[
H^*(B\G_2(q); \mathbb{F}_2)\cong \mathbb{F}_2[d_4,d_6,d_7,y_3,y_5]/I
\]
where $I$ is the ideal generated by the polynomials
\begin{enumerate}
 \item[] $y_5^2+y_3d_7+y_3^2d_4$
 \item[] $y_3^4+y_5d_7+y_3^2d_6$.
\end{enumerate}
The action of $\A_2$ is determined by
\[
Sq^2d_4=d_6,\quad Sq^1d_6=d_7,\quad Sq^1y_3=0,\quad
Sq^2y_3=y_5,\quad Sq^1y_5=u_3^2.
\]
\end{theorem}
To specify a finite (untwisted) group of Lie type $\G(q)$ we need to
know a compact connected Lie group $\G$ and a finite field $\fq$.
The following theorem of Quillen~\cite{Qu} and Friedlander~\cite{Fr}
relates the classifying space of the finite group $\G(q)$ to the
classifying space of $\G$.
\begin{theorem}[Quillen, Friedlander]
\label{flie} For every prime $l$ not dividing $q$ there is a
homotopy pullback diagram
\begin{equation}
\label{dgrmflt} \xymatrix{
B\G(q)^{\wedge}_l\ar[r]\ar[d] & B\G^{\wedge}_l\ar[d]^{\Delta}\\
B\G^{\wedge}_l \ar[r]^-{1\times\Psi^q} & B\G^{\wedge}_l\times
B\G^{\wedge}_l.}
\end{equation}
\end{theorem}
Consider the fibration sequence $(\Omega B\G_2)^{\wedge}_2\lra
(B\G_2(q))^{\wedge}_2\lra (B\G_2)^{\wedge}_2$. The cohomology
$H^*(B\G_2)$ is the mod~$2$ Dickson algebra of rank~$3$, while the
cohomology of the fibre at the prime two is $P[u_3]/u_3^4\otimes
E[u_5]$ where $u_5=Sq^2 u_3$ (see for example \cite[Theorem
6.2]{MT}). Applying Smith's Big Collapse Theorem \cite{Sm} to the
Eilenberg--Moore spectral sequence calculation, we have
\[
H^*(B\G_2(q))\cong P[d_4,d_6,d_7] \otimes P[u_3]/u_3^4\otimes E[u_5]
\]
as $P[d_4,d_6,d_7]$--modules, with the action of $\A_2$ given by
$Sq^2d_4=d_6$, ${Sq^1d_6=d_7}$, $Sq^2u_3=u_5$, $Sq^1u_5=u_3^2$.

\begin{lemma}
The natural map $H^*(B\G_2)\lra H^*(B\G_2(q))$ has a left inverse in
the category of unstable algebras over $\A_2$.
\end{lemma}
\begin{proof}
To show that the natural map $H^*(B\G_2)\lra H^*(B\G_2(q))$ has a
left inverse in the category of unstable algebras over $\A_2$ we
proceed as in the proof of Proposition \ref{Lannes}.

There is a bijection \cite{La}
\[
[B(\mathbb{Z}/2)^3, B\G_2]\cong\Hom_{\A_2}(H^*(B\G_2),
H^*(B(\mathbb{Z}/2)^3))
\]
which assigns to the inclusion of the mod~$2$ Dickson invariants of
rank~$3$ into $H^*((\mathbb{Z}/2)^3)$ a geometric map $f\colon
B(\mathbb{Z}/2)^3\lra B\G_2$. Now proceed exactly as in the proof of
Proposition \ref{Lannes} using the pullback diagram
\[
\xymatrix{B\G_2(q)\ar[d]^(.4){p}\ar[r] & B\G_2\ar[d]^{\Delta}\\
B\G_2\ar[r]^-{(1,\Psi^q)} & B\G_2\times B\G_2}
\]
in place of pullback diagram \eqref{maindgrm}.
\end{proof}

As we shall show, there are two possible algebra extensions
\[
\mathbb{F}_2\lra H^*(B\G_2)\lra H^*(B\G_2(q))\lra H^*(\Omega
B\G_2)\lra \mathbb{F}_2
\]
determined by the value of a single parameter $A\in \mathbb{F}_2$.

Let us consider possible multiplications on $H^*(B\G_2(q))$. We
denote by $\DI(3)$ the loop space $\Omega B\G_2\simeq\G_2$. Denote
by $y_3\in H^*(B\G_2(q))$ the unique generator which maps to $u_3\in
H^*(\DI(3))$ under the map induced by the fibre inclusion. Let
$y_5\in H^*(B\G_2(q))$ denote $Sq^2(y_3)$. The algebra structure of
$H^*(B\G_2(q))$ splits in the same way as a module structure if and
only if the subalgebra of $H^*(B\G_2(q))$ generated by $y_3$ and
$y_5$ is isomorphic to $P[u_3]/(u_3^4)\otimes E[u_5]$.

Looking at the action of the Steenrod algebra $\A_2$, we have $
y_5^2 = Sq^4(y_3^2)$. Because of dimensional reasons, $y_5^2$ can be
presented as
\begin{equation}
\label{rel2} y_5^2 = Au_3d_7 + By_3^2d_4 + Cd_4d_6.
\end{equation}
On the one hand we have $Sq^1 y_5^2=0$, and on the other hand
applying $Sq^1$ to relation \eqref{rel2}, we have $Sq^1y_5^2 =
Cd_4d_6$. This implies that $C=0$. Next
\[
y_3^4 = Sq^2y_5^2 = Ay_5d_7 + By_3^2d_6
\]
and
\[
0=Sq^1y_3^4 = (A+B)y_3^2d_7.
\]
This implies that $A=B$. Therefore all relations in cohomology of
$\G_2(q)$ depend on a single parameter $A\in\mathbb{F}_2$. Notice
that if $A=0$, then there is a ring homomorphism $H^*(\DI(3))\lra
H^*(B\G_2(q))$ and the composite
\begin{equation}
\label{splitting2} H^*(B\DI(3))\otimes H^*(\DI(3))\lra
H^*(B\G_2(q))^{\otimes 2} \lra H^*(B\G_2(q))
\end{equation}
is an isomorphism of vector spaces and a map of algebras. Hence
composite \eqref{splitting2} is an isomorphism of algebras over the
Steenrod algebra. The same argument as in the case of
$H^*(B\Sol(q))$ in Lemma \ref{product} shows that the cohomology
algebra cannot be split as the tensor product of $H^*(B\G_2)$ and
$H^*(\Omega B\G_2)$ since $H^*(\DI(3))\ncong H^*((\mathbb{Z}/2)^3)$.
Therefore $A$ must be $1$. This proves the algebra structure of the
mod~$2$ cohomology of $B\G_2(q)$ stated in Theorem~\ref{cohg2q}.

\end{document}